# Considerations on M/G/∞ Queue Busy Period Variance

**Manuel Alberto M. Ferreira**
Instituto Universitário de Lisboa (ISCTE – IUL), ISTAR- IUL
Lisboa, Portugal
manuel.ferreira@iscte.pt

**ABSTRACT**

Some important results on the variance of the M/G/∞ queue busy period are presented. Often, this parameter depends on the whole structure of the service time distribution. So, the importance of the bounds presented, depending only on some parameters. Also, some bounds for service time distributions important in reliability theory, with technological and financial applications, are presented and a particular stochastic ordination is explored to explain those results. Finally, an interesting result about the exponentiality of the M/G/∞ queue busy period, resulting from the coefficient of variation study is shown.

## 1. BOUNDS FOR THE VARIANCE OF THE M/G/∞ QUEUE BUSY PERIOD

The $M/G/\infty$ busy period variance is given by, see for instance (1),

$$VAR[B] = \frac{2e^\rho \int_0^\infty (e^{\lambda \int_t^\infty [1-G(v)]dv} - 1)dt}{\lambda} - \left(\frac{e^\rho - 1}{\lambda}\right)^2 \quad (1.1)$$

where $\lambda$ is the customers arrival rate, $\rho$ is the traffic intensity, given by

$$\rho = \lambda \alpha \quad (1.2),$$

being α mean service time, and G(.) is the service time distribution function. $B$ designates the busy period random variable. So, in general, $VAR[B]$ depends on the whole structure of the service time distribution.

In (2) are deduced lower and upper bounds for the $M/G/\infty$ busy period variance that depend only on $\rho$, $\lambda$ and $\gamma_s$ - the service time coefficient of variation – that are then improved in (3) assuming the form

$$\lambda^{-2} \max[e^{2\rho} + e^\rho \rho^2 \gamma_s^2 - 2\rho e^\rho - 1, 0] \leq VAR[B] \leq \lambda^{-2}\{2e^\rho(\gamma_s^2 + 1)(e^\rho - 1 - \rho) - (e^\rho - 1)^2\} \quad (1.3).$$

Note that

- With $\gamma_s = 0$ the bounds are coincident and are $VAR[B]$ value for the $M/D/\infty$ queue system – constant services with value α,

- For the $M/M/\infty$ queue system – exponential service times

$$VAR[B^M] = \frac{2e^\rho \left(1 + \rho \sum_{n=1}^{\infty} \frac{\rho^n}{nn!}\right) - e^{2\rho} - 1}{\lambda^2} \qquad (1.4).$$

In (3) it is shown that (1.1) is equivalent to

$$VAR[B] = \frac{e^{2\rho} - 2\rho e^\rho - 1}{\lambda^2} + e^\rho \sigma_s^2 + \frac{e^\rho}{\lambda^2} \sum_{n=3}^{\infty} \frac{\rho^n}{n!}\left(b_{n-2} - \frac{2}{1+\gamma_s^2}\right) \qquad (1.5)$$

where $\sigma_s^2$ and $\gamma_s$ are, respectively, the variance and the coefficient of variation of the service time and

$$b_n = \frac{2(n+2) \int_0^\infty [1 - G(v)]^{n+1} dv}{\alpha^{n+2}(1+\gamma_s^2)}, n = 0,1,2,\ldots \qquad (1.6).$$

As $\frac{2}{1+\gamma_s^2} \leq b_n \leq 2, n = 0,1,2,\ldots$ then follow the bounds in formula (1.3). But, for the $M/M/\infty$ system, the upper bound may be improved:

- $b_n = \dfrac{(n+2) \int_0^\infty \left[\int_v^\infty e^{-\frac{x}{\alpha}} dx\right]^{n+1} dv}{\alpha^{n+2}} = \dfrac{(n+2) \int_0^\infty \left(\left[-\alpha e^{-\frac{x}{\alpha}}\right]_v^\infty\right)^{n+1} dv}{\alpha^{n+2}} =$

$\dfrac{\left(n+2 \int_0^\infty \left(\alpha e^{-\frac{v}{\alpha}}\right)^{n+1} dv\right)}{\alpha^{n+2}} = (n+2) \dfrac{\alpha^{n+1}}{\alpha^{n+2}} \int_0^\infty e^{-\frac{(n+1)v}{\alpha}} dv = \dfrac{n+1}{n+2}.$

- So,

$$VAR[B] = \frac{e^{2\rho} - 2\rho e^\rho - 1}{\lambda^2} + e^\rho \alpha^2 + \frac{e^\rho}{\lambda^2} \sum_{n=3}^{\infty} \frac{\rho^n}{n!}\left(\frac{n}{n-1} - 1\right) =$$
$$\frac{e^{2\rho} - 2\rho e^\rho - 1}{\lambda^2} + e^\rho \alpha^2 + \frac{e^\rho}{\lambda^2} \sum_{n=3}^{\infty} \frac{\rho^n}{n!} \frac{1}{n-1}.$$

- But

$$\sum_{n=3}^{\infty} \frac{\rho^n}{n!} \frac{1}{n-1} \leq \sum_{n=3}^{M} \frac{\rho^n}{n!} \frac{1}{n-1} + \frac{1}{M} \sum_{n=M+1}^{\infty} \frac{\rho^n}{n!} = \sum_{n=3}^{M} \frac{\rho^n}{n!} \frac{1}{n-1} + \frac{1}{M}\left(e^{\rho} - 1 - \rho - \cdots - \frac{\rho^M}{M!}\right).$$

- And

$$VAR[B^M] \leq \frac{e^{2\rho} - 2\rho e^{\rho} - 1}{\lambda^2} + e^{\rho}\alpha^2$$
$$+ \frac{e^{\rho}}{\lambda^2}\left(\sum_{n=3}^{M} \frac{\rho^n}{n!} \frac{1}{n-1} + \frac{1}{M}\left(e^{\rho} - \sum_{n=0}^{M} \frac{\rho^n}{n!}\right)\right) \quad (1.7).$$

The upper bound given by (1.7) is better than the one of (1.3) if a $M$ great enough is considered. Computing (1.3) and (1.7), see Table 1, for $\lambda = 1$ and $\rho = 0.5; 1; 10; 20; 50; 100$, for instance, it is noted that for $M \geq 14$ the upper bound given by (1.7) is better than the one given by (1.3):

**Table 1.** Computation of upper and lower bounds given by (1.3) and (1.7), for the $M/M\infty$ system with $\lambda = 1$.

| $\rho$ | Upper Bound (1.3) | Upper Bound (1.7); $M = 14$ | Lower Bound (1.3) |
|---|---|---|---|
| 0.5 | 0.55954328 | 0.50046123 | 0.48174095 |
| 1 | 4.8574775 | 3.9415704 | 3.6707743 |
| 10 | $1.4545705 \times 10^9$ | $5.6362048 \times 10^8$ | $4.8692729 \times 10^8$ |
| 20 | $7.061558 \times 10^{17}$ | $2.5325047 \times 10^{17}$ | $2.3538545 \times 10^{17}$ |
| 50 | $8.0643512 \times 10^{43}$ | $2.8801252 \times 10^{43}$ | $2.6881171 \times 10^{43}$ |
| 100 | $2.167792 \times 10^{87}$ | $7.7421139 \times 10^{86}$ | $7.2259735 \times 10^{86}$ |

For service time distributions related with the exponential distribution, important in reliability theory, see (4), it is easy to show that:

- Service time distribution NBUE-New Better than Used in Expectation with mean α

$$VAR[B^{NBUE}] \leq VAR[B^M] \quad (1.8),$$

- Service time distribution NWUE-New Worse than Used in Expectation with mean α

$$VAR[B^{NWUE}] \geq VAR[B^M] \quad (1.9),$$

- Service time distribution DFR-Decreasing Failure Rate

$$VAR[B^{DFR}] \geq \frac{2e^{\rho}\rho \sum_{n=1}^{\infty} \frac{\rho^n}{nn!} e^{-\frac{n}{2}(\gamma_s^2-1)} + e^{\rho}(2 - e^{\rho} - e^{-\rho})}{\lambda^2} \qquad (1.10),$$

- Service time distribution IMRL- Increasing Mean Residual Life

$$VAR[B^{IMRL}] \geq \frac{e^{\rho}}{\rho}\mu_2 \sum_{n=1}^{\infty} \frac{\rho^n}{nn!} e^{-n\left(\frac{2\alpha\mu_3}{3\mu_2^2}-1\right)} + \lambda^{-2} e^{\rho}(2 - e^{\rho} - e^{-\rho}) \qquad (1.11)$$

where $\mu_r$ are the service time distribution order $r$ moments centered in the origin, $r=2,3$.

Note that, in this context of reliability theory, the service time may be the lifetime of a technological device or of a financial asset, for instance.

## 2. VARIABILITY ORDINATION

The variability ordination between two random variables, $G_1$ and $G_2$, is symbolized by

$$G_1 \underset{V}{\leq} G_2 \qquad (2.1)$$

and it means that

$$\begin{aligned} &a) \ E[G_1] = E[G_2] \\ &b) \int_t^{\infty} [1 - G_1(v)]dv \leq \int_t^{\infty} [1 - G_2(v)]dv, t \in \mathbb{R} \end{aligned} \qquad (2.2),$$

see (5).

It is possible to show that
- If in $M/G_1/\infty$ and $M/G_2/\infty$,

$$G_1 \underset{V}{\leq} G_2$$

so

$$VAR[B^{G_1}] \leq VAR[B^{G_2}]$$

and, as a conjecture, that

$$B_1 \underset{V}{\leq} B_2$$

being $B_1$ and $B_2$ the busy period random variables associated to $M/G_1/\infty$ and $M/G_2/\infty$ respectively.

This result has the formulae (1.8) and (1.9) as particular cases.

## 3. COEFFICIENT OF VARIATION OF THE M/G/∞ QUEUE BUSY PERIOD

After the results obtained for the variance, result evidently similar results for the coefficient of variation. Not being interesting to present a sum of trivial results, only one non-trivial result will be presented in this section. Call $\gamma_B$ the coefficient of variation of the $M/G/\infty$ queue busy period. It fulfills the equality

$$\gamma_B^2 = \frac{2e^\rho \lambda \beta}{(e^\rho - 1)^2} - 1, \text{with } \beta = e^\rho \int_0^\infty \left[ e^{-\lambda \int_0^t [1-G(v)]dv} - e^{-\rho} \right] dt \quad (3.1).$$

With given $\lambda$ and α great enough, for service time distributions such that $1 - G(t) \cong 1$, $\beta \cong e^\rho \int_0^\infty e^{-\lambda t} dt = \frac{e^\rho}{\lambda}$ and so $\gamma_B^2 \cong 1$, indicating exponential behavior for the busy period distribution in these conditions, see (1).

## 4. CONCLUSIONS

The exigency of infinite servers means, in general, that either there is no distinction between the customer and server or when a customer arrives meets immediately an available server. It is not imposed the physical presence of infinite servers. So, in some situations, the manager of the queue system must take in account that he/she must have a certain reserve of servers to join the service system. In consequence the variance of the length of the M/G/∞ queue busy period is an important parameter to consider. Great variances of course imply a great disturbance in the costs of the service because they may imply an over dimension of the reserve of the servers.